\documentclass[12pt]{amsart}

\usepackage{a4wide}
\usepackage{amssymb}
\usepackage{hyperref}

\theoremstyle{plain}
\newtheorem{Theorem}{Theorem}

\newtheorem{Proposition}[Theorem]{Proposition}

\newtheorem{Question}[Theorem]{Question}
\newtheorem{Conjecture}[Theorem]{Conjecture}
\newtheorem{ourque}{Question}

\newcommand{\Z}{{\mathbb{Z}}}
\newcommand{\Q}{{\mathbb{Q}}}
\newcommand{\Gal}{{\rm Gal}}
\newcommand{\bfsigma}{\mbox{\boldmath$\sigma$}}

\def\hefresh{
   \def\hefreshD{\mathop{\raise1.5pt\hbox{${\smallsetminus}$}}}
   \def\hefreshS{\mathop{\raise0.85pt\hbox{$\scriptstyle\smallsetminus$}}}
   \mathchoice{\hefreshD}{\hefreshD}{\hefreshS}{\hefreshS}}
\renewcommand{\setminus}{\hefresh}

\numberwithin{Theorem}{section}

\begin{document}

\bibliographystyle{alpha}

\title[Open Problems in the Theory of Ample Fields]{Open Problems in the Theory of Ample Fields}
\author{Lior Bary-Soroker}
\address{Institut f\"ur Experimentelle Mathematik, Universit\"at Duisburg-Essen, Ellernstr.~29, 45326 Essen, Germany}
\email{lior.bary-soroker@uni-due.de}
\author{Arno Fehm}
\address{Universit\"at Konstanz,
Fachbereich Mathematik und Statistik,
Fach D 203,
78457 Konstanz,
Germany}
\email{arno.fehm@uni-konstanz.de}

\begin{abstract}
Fifteen years after their discovery, ample fields now stand at the center of research in contemporary Galois theory and attract more and more attention also from other areas of mathematics. This survey gives an introduction to the theory of ample fields and discusses open problems.
\end{abstract}

\maketitle

\section{Introduction}

\noindent
In the middle of the 1990's, Pop realized that all fields
for which a certain Galois theoretic conjecture was proven 
(namely the regular solvability of finite split embedding problems)
share a common property: the set of rational points of any smooth curve over such a field is either empty of infinite.
Moreover, in \cite{Pop} he showed that the conjecture actually holds for all fields satisfying this property.
Since then these fields, nowadays called {\em ample} fields\footnote{Some authors prefer the term {\em large} field.} -- the topic of this survey -- play a central role in Galois theory.

There are several equivalent ways to define ample fields,
and this notion captures well the intuitive concept of a `large' field.
For example, the class of ample fields subsumes several seemingly unrelated classes of large fields,
like algebraically closed, real closed, separably closed, and Henselian valued fields.
A typical example of an ample field that does not fall into any of these categories is the field of totally $S$-adic numbers $\Q_{{\rm tot},S}$, $S$ a finite set of places of $\Q$ --
the maximal extension of $\Q$ in which each place in $S$ is totally split.
For example, the field of totally real algebraic numbers, which plays an important role in number theory, is of this kind.

In recent years, ample fields 
attracted more and more attention also in other subjects,
like arithmetic geometry \cite{Debes,Kollar,MoretBailly2001,FehmPetersen}, 
valuation theory \cite{KuhlmannPlaces,AKP}, 
and model theory \cite{Koenigsmann,Tressl,JunkerKoenigsmann}.

In this note we survey the basics of the theory of ample fields and discuss open problems.
The Galois-theoretic aspects of the theory of ample fields are covered well in the literature,
and we refer the reader to the beautiful survey paper \cite{DD} of D\`ebes and Deschamps,
and the comprehensive book \cite{JardenAP} by Jarden.
We will focus on new developments that took place since the publication of \cite{DD},
and on connections between the theory of ample fields and other subjects.

\section{Background}

\subsection{Characterization}
Apart from being given different names (like large, anti-Mordellic, fertile, pop)
there are also several equivalent ways to define 
the class of ample fields, cf.~\cite{Pop}:

\begin{Proposition}
The following properties for a field $K$ are equivalent:
\begin{enumerate}
\item[(1)] Every smooth $K$-curve with a $K$-rational point has infinitely many such points.
\item[(2)] If $V$ is a smooth $K$-variety, then the set of $K$-rational
   points of $V$ is either empty or Zariski-dense in $V$.
\item[(3)] $K$ is existentially closed in the field of formal Laurent series $K((t))$.
\item[(4)] $K$ is existentially closed in some field extension that admits a nontrivial Henselian valuation.
\end{enumerate}
\end{Proposition}

Here, $K$ is called {\em existentially closed} in an extension $F$ if
every existential first-order sentence in the language of rings with parameters from $K$ which holds in $F$, also holds in $K$.
Equivalently, $K$ is existentially closed in $F$ if for any $x_1,\dots,x_n\in F$ there exists a $K$-homomorphism from the ring $K[x_1,\dots,x_n]$ to $K$.
An {\bf ample} field is a field $K$ that satisfies the equivalent conditions (1)-(4).

\subsection{Ample and non-ample fields}

The most important known properties of fields that imply ampleness can be grouped into three classes:

\begin{enumerate}
\item arithmetic-geometric
\item topological
\item Galois-theoretic
\end{enumerate}

{\em Arithmetic-geometric properties:} Every algebraically closed, separably closed, and more generally, pseudo algebraically closed field, is ample. 
Here, a field $K$ is called pseudo algebraically closed (PAC, see \cite{FJ})
if every absolutely irreducible $K$-variety has a $K$-rational point.
Even more generally, fields that satisfy a geometric local-global principle for varieties are ample, \cite{Pop}.
An example of such fields are the ${\rm P}S{\rm C}$ fields:
If $S$ is a finite set of places of $\Q$, then a field $K\subseteq\Q_{{\rm tot},S}$ is called ${\rm P}S{\rm C}$
if every smooth $\Q$-variety that has a point over $\Q_p$ for each $p\in S$ (with $\Q_\infty=\mathbb{R}$)
has a $K$-rational point.

{\em Topological properties:} Every field which is complete with respect to a nontrivial absolute value is ample. Also, every field which admits a nontrivial Henselian valuation is ample.
This can be generalized further to the quotient fields of domains that are complete (more generally Henselian) with respect to an ideal \cite{PopHenselian}, or a norm \cite{FehmParan}.

{\em Galois-theoretic properties} that imply ampleness are discussed in Section \ref{sec:Galois}.

On the other side, there are three basic classes of non-ample fields:

\begin{enumerate}
\item finite fields
\item number fields, i.e.~finite extensions of $\Q$
\item function fields, i.e.~fields that are a finitely generated and transcendental extension of another field
\end{enumerate}

In particular, all global fields are non-ample, in contrast to the fact that all local fields are ample.
Apart from these three classes, only very few non-ample fields are known.
As D\`ebes puts it in \cite{Debes}:

\begin{quote}
\textquotedblleft [...] it happens to be difficult to produce non-ample fields at all
(at least inside $\bar{\mathbb{Q}}$ and apart from number fields)\textquotedblright.
\end{quote}

For some examples of non-ample fields that do not fall into any of the above three classes
see \cite{KoenigsmannIGP,LozanoRobledo,AFembeddings}.
Surprisingly, to the best of our knowledge,
all infinite non-ample fields appearing in the literature are Hilbertian.

\subsection{Properties}

As explained before,
the notion of ample fields captures very well the intuitive concept of a large field.
One of the properties every notion of large fields should certainly satisfy 
is the following, cf.~\cite{Pop}:

\begin{Proposition}
The class of ample fields is closed under algebraic extensions.
\end{Proposition}

In addition we have:

\begin{Proposition}
The class of ample fields is an elementary class. In particular,
it is closed under elementary equivalence in the language of rings.
\end{Proposition}

\subsection{Abelian varieties}

Let $A$ be an abelian variety defined over a finitely generated field $K$.
It is known that the Mordell-Weil group $A(K)$ is finitely generated,
and in particular 
$$
 {\rm rank}(A(K)):={\rm dim}_{\mathbb{Q}}(A(K)\otimes\mathbb{Q}),
$$
the {\bf rank} of $A$ over $K$, is finite.
On the other hand, if $K$ is algebraically closed and not algebraic over a finite field, then ${\rm rank}(A(K))=\infty$.
At least in characteristic zero, this holds more generally for arbitrary ample fields, cf.~\cite{FehmPetersen}:

\begin{Proposition}\label{FehmPetersen}
Let $K$ be an ample field of characteristic zero and $A/K$ a non-zero abelian variety.
Then ${\rm rank}(A(K))=\infty$.
\end{Proposition}

Although one can construct a non-ample field
over which every abe\-li\-an variety has infinite rank
(i.e.~the converse of Proposition \ref{FehmPetersen} does not hold),
we do not know any natural example of this kind,
and none which is algebraic over $\mathbb{Q}$.

\section{Algebraic fields}

\noindent
We denote by  $\mathbb{Q}_{\rm ab}$ the maximal abelian extension of $\mathbb{Q}$,
which, by the Kronecker-Weber theorem, coincides with the maximal cyclotomic extension
$\Q(\zeta_n:n\in\mathbb{N})$ of $\Q$.

\subsection{}
Our first question can be simply stated as follows. 

\begin{ourque}\label{ques:Qab}
Is $\mathbb{Q}_{\rm ab}$ ample? 
\end{ourque}

The importance of this question lies in the fact that, via Pop's results on ample fields,
a positive answer would give a proof of the following conjecture of Shafarevich, cf.~\cite{DD}:

\begin{Conjecture}[Shafarevich]\label{con:Sha}
The absolute Galois group $G_{\mathbb{Q}_{\rm ab}}$ of  $\mathbb{Q}_{\rm ab}$ is a free profinite group.
\end{Conjecture}

As evidence for this conjecture
one has Iwasawa's theorem 
that the maximal prosolvable quotient of  $G_{\mathbb{Q}_{\rm ab}}$ is prosolvable free,
and Tate's result that the cohomological dimension of $G_{\mathbb{Q}_{\rm ab}}$ is $1$. 
Another piece of evidence is that the geometric Shafarevich conjecture --
the analogue of the Shafarevich conjecture for global function fields -- holds true, \cite{Harbater,Popetale}.
A proof of Conjecture \ref{con:Sha} could be seen as a step towards understanding the absolute Galois group of $\Q$,
since $G_\Q$ is an extension of the well-known group $\Gal(\Q_{\rm ab}/\Q)\cong\hat{\mathbb{Z}}^\times$ by $G_{\Q_{\rm ab}}$.

\subsection{}\label{sec:Qsolv}
It is interesting to note that it is even unknown whether the bigger field $\Q_{\rm solv}$,
the maximal solvable Galois extension of $\Q$, is ample.
Even worse, while Frey proved that $\Q_{\rm ab}$ is {\em not} PAC, \cite[Corollary 11.5.7]{FJ},
it is a long-standing open question whether $\Q_{\rm solv}$ is PAC,
equivalently, whether every absolutely irreducible $\Q$-variety admits a solvable rational point.

The ampleness of $\Q_{{\rm tot},S}$ mentioned in the introduction
follows from the fact that $\Q_{{\rm tot},S}$ is ${\rm P}S{\rm C}$, \cite{MoretBailly,GPR}.
This fact plays an important role in the study of potential modularity of Galois representations,
which motivated Taylor to ask in \cite{Taylor} whether also 
$\Q_{{\rm tot},S}\cap\Q_{\rm solv}$ is ${\rm P}S{\rm C}$. 
A positive answer to this would have \textquotedblleft extremely important consequences\textquotedblright, as he points out,
and of course it would imply that $\Q_{{\rm tot},S}\cap\Q_{\rm solv}$ and consequently also $\Q_{\rm solv}$
are ample.

\subsection{}
A positive answer to Question~\ref{ques:Qab} would also settle the following problem.
In their 1974 paper \cite{FreyJarden}, Frey and Jarden
investigate the rank of elliptic curves over certain fields which are not finitely generated.
For example, they prove that any elliptic curve $E/\mathbb{Q}$
acquires infinite rank over the field $\mathbb{Q}(\sqrt{n}:n\in\mathbb{Z})$,
so in particular
${\rm rank}(E(\mathbb{Q}_{\rm ab}))=\infty$.
This led them to ask for a generalization of this result to arbitrary abelian varieties:

\begin{Question}[Frey-Jarden]\label{FreyJarden}
Is ${\rm rank}(A(\mathbb{Q}_{\rm ab}))=\infty$ for every non-zero abelian variety $A$ over $\mathbb{Q}$?
\end{Question}

However, this generalization turned out to be complicated and created lots of activity leading
to positive partial solutions for various kinds of abelian varieties, e.g.~Jacobians of curves
with certain properties, see for example \cite{RosenWong,Petersen} and the references therein.
Nevertheless, if $\mathbb{Q}_{\rm ab}$ is an ample field, then Proposition \ref{FehmPetersen}
would immediately imply that the answer to Question~\ref{FreyJarden} is positive for all abelian varieties.

\section{Galois theory of ample fields}\label{sec:Galois}

\noindent
The main topic of {\em Field Arithmetic} is the interplay of Galois theoretic and arithmetic
properties of fields.
An example for this interplay is the classical Artin-Schreier theorem, which says that
the absolute Galois group $G_K$ of a field $K$ is finite if and only if
$K$ is separably closed (in which case $G_K=1$) or real closed (in which case $G_K\cong\Z/2\Z$).
Therefore, if a field has a finite absolute Galois group, then it is ample.
This section deals with other conditions on $G_K$ that imply that $K$ is ample.

Note that there cannot be any conditions on $G_K$ implying that $K$ is {\em non-ample},
since every profinite group that occurs as the absolute Galois group of a field
also occurs as the absolute Galois group of an ample field:
Every absolute Galois group is also an absolute Galois group of a field $K$ of characteristic zero (\cite[Corollary 22.2.3]{Efrat}),
and $G_K$ is isomorphic to a closed subgroup of the absolute Galois group of $K((t))$;
therefore, the ample field $K((t))$ has an algebraic, and consequently also ample extension $K'$ with $G_K\cong G_{K'}$.
\subsection{}

Colliot-Th\'el\`ene observed that a field $K$ is ample whenever $G_K$ is a pro-$p$ group for some
prime number $p$, \cite{ColliotThelene,JardenAmple}.
This covers again in particular
the archimedean local fields $\mathbb{R}$ and $\mathbb{C}$.
Since pro-$p$ groups are prosolvable, and all local fields have prosolvable absolute Galois groups
we want to ask for the following generalization of Colliot-Th\'el\`ene's result:

\begin{ourque}\label{ques:prosolv}
Is every infinite field $K$ with prosolvable absolute Galois group $G_K$ ample?
\end{ourque}

\subsection{}\label{sec:fg}

Another Galois theoretic property that the local fields $\Q_p$, $\mathbb{R}$, and $\mathbb{C}$
share, apart from having a prosovable absolute Galois group,
is that their absolute Galois groups are all finitely generated (as topological groups).
The following question occurs as a conjecture in \cite{JunkerKoenigsmann}:
\begin{ourque}\label{ques:fg}
Is every infinite field $K$ with finitely generated absolute Galois group $G_K$ ample?
\end{ourque}

As evidence towards a positive answer let us recall the following result of Jarden:
If $F$ is a countable Hilbertian field and $e\in\mathbb{N}$,
then the subset $\Sigma$ of those $\bfsigma\in (G_{F})^e$ whose fixed field $F_{\rm sep}(\bfsigma)$
is PAC has Haar measure $1$. Question~\ref{ques:fg} is therefore equivalent to asking
whether $K=F_{\rm sep}(\bfsigma)$ is ample also for the $\bfsigma$ in the complement of $\Sigma$ --
which has measure zero. Note that this complement contains all $\bfsigma$
for which the fixed field is real closed, and all for which it is Henselian but not separably closed.

Other pieces of evidence for a positive answer to both Question~\ref{ques:prosolv} and Question~\ref{ques:fg}
are known conditions on $G_K$ that imply that $K$ is Henselian and thus ample.
For example, if $G_K$ is solvable and $G_{K(\sqrt{-1})}$ is not projective, 
then $K$ admits a nontrivial Henselian valuation,
\cite{Koenigsmannsolvable}.
Further evidence will be given in the following sections.

\subsection{}

By Proposition~\ref{FehmPetersen}, a positive answer to Question~\ref{ques:fg} would also settle the following conjecture 
from \cite{Larsen}:
\begin{Conjecture}[Larsen]\label{Larsen}
Let $K/\Q$ be an algebraic extension with $G_K$ finitely generated,
and let $A$ be a non-zero abelian variety over $K$. Then ${\rm rank}(A(K)) = \infty$.
\end{Conjecture}
It is interesting to note that certain cases of this conjecture are already proven,
for example the case 
that $G_K$ is procyclic,
see \cite{ImLarsen}.
This can in turn be seen as additional evidence for a positive answer of Question~\ref{ques:fg}.

\section{Virtually ample fields}

\noindent
A field $K$ is {\bf virtually ample} if some finite extension of $K$ is ample.
Exchanging quantifiers in this definition we get a weaker condition:
A field $K$ is called {\bf weakly ample} if for every $K$-curve $C$
there exists a finite extension $L$ of $K$ such that $C(L)$ is infinite.

Trivially, every ample field is virtually ample, and every virtually ample field is weakly ample.
One could think that every virtually ample field should be ample,
but no proof is known. 
As evidence one usually refers to the fact that similar descent properties hold for
Henselian valued fields.
We want to ask more generally:

\begin{ourque}\label{ques:weakly}
What is the relation between ample, virtually ample, and weakly ample fields?
Is every weakly ample field virtually ample? Is every virtually ample field ample?
\end{ourque}


\subsection{}
If every virtually ample field is ample, then one would immediately get that every 
number field is non-ample, since $\Q$ is. The fact that number fields are non-ample is of course
well-known, but all the proofs known to us involve deep arithmetic-geometric arguments.
On the other hand, already Fermat knew that for example the equation $x^4+y^4=z^4$ has no nontrivial integer solutions,
which gives an easy proof that $\Q$ is non-ample.

\subsection{}
An important class of weakly ample fields are the fields considered in Section \ref{sec:fg}:
Infinite fields with finitely generated absolute Galois group.
Indeed, a $K$-curve of gonality $\gamma$ has infinitely many points of residue degree bounded by $\gamma$,
but if $G_K$ is finitely generated then $K$ has only finitely many extensions of degree $\gamma$,
so there must be one over which infinitely many points become rational.
Thus if every weakly ample field is ample, Question~\ref{ques:fg} would have a positive answer.

\section{Radically closed fields}

\noindent
A field $K$ is called {\bf radically closed} if the multiplicative group $K^\times$ is divisible.
Any radically closed field contains all roots of unity and hence, by Kummer theory, has no solvable finite Galois extensions.
One would like to know the following:

\begin{ourque}\label{ques:radically}
Is every radically closed field ample?
\end{ourque}

\subsection{}
For example, $\mathbb{Q}_{\rm solv}$ is radically closed,
so a positive answer to Question~\ref{ques:radically} would answer the question
mentioned in Section \ref{sec:Qsolv}, wheth\-er $\mathbb{Q}_{\rm solv}$ is ample. 
It would also prove a long-standing conjecture in model theoretic algebra:

\subsection{}
An infinite field $K$ is called {\bf minimal} if every subset $A\subseteq K$ which is definable
by a first-order formula in the language of rings is finite or cofinite.
For example, algebraically closed fields have this property,
and an old conjecture of Podewski says that they should be the only ones:
\begin{Conjecture}[Podewski]
Every minimal field is algebraically closed.
\end{Conjecture}
Surprisingly, so far this is proven only for minimal fields of {\em positive} characteristic, \cite{Wagner}.
As for the general case, it seems that it is even not known whether there exists an infinite field $K$
and a polynomial $f\in K[X]$ such that the definable set $K\setminus f(K)$ is non-empty but finite.

Koenigsmann observed that a positive answer to Podewski's conjecture would follow
from a positive answer to Question~\ref{ques:radically}.
Indeed, minimal fields are easily seen to be radically closed,
and an ample and minimal field $K$ must already be algebraically closed:
If $f\in K[X]$ is an irreducible and separable polynomial, then the image $f(K)$ is cofinite in $K$.
However, Krasner's lemma implies that in the Henselian field $K((t))$ there exist infinitely many elements $x$
for which $f(X)-x$ generates the same extension as $f(X)$.
Since by the existential closedness of $K$ in $K((t))$ the same holds in $K$
one gets a contradiction to the assumption that $f(K)$ is cofinite.

Note that since a Hilbertian field is never radically closed,
this also implies that a minimal field which is not algebraically closed would have to be both
non-ample and non-Hilbertian -- as mentioned above we do not have any example of a field that satisfies these two properties simultaneously.

\subsection{}
Call $K$ {\bf almost radically closed} if for every $n$
the subgroup of $n$-th powers has finite index in $K^\times$.
One could strengthen the above question and ask whether every infinite almost radically closed field is ample.
Since every field with a finitely generated absolute Galois group is
almost radically closed, this would give an answer to Question~\ref{ques:fg}.

\section{The conjecture of D\`ebes and Deschamps}

\noindent
We want to conclude with some explanations concerning 
the Galois theoretic conjecture  mentioned in the introduction.
This conjecture is a weakening of a conjecture of D\`ebes and Deschamps that was first published in \cite{DD}.

The classical {\bf inverse Galois problem}, which dates back to the 19th century, asks whether every finite group occurs as a Galois group over $\Q$.
By Hilbert's irreducibility theorem, this is equivalent to asking the same question for the rational function field $\Q(t)$.
Consider the following two generalizations of the inverse Galois problem:

The {\bf regular inverse Galois problem} for a field $K$ asks whether for every finite group 
$G$ there exists a Galois extension $F$ of $K(t)$ with Galois group $G$ such that $F/K$ is regular.
It was conjectured that this holds true for any field $K$. 

A {\bf finite split embedding problem} over a field $E$ is a finite Galois extension $F_1/E$ together with
an epimorphism of finite groups
$\alpha: G\rightarrow\Gal(F_1/E)$
that splits, i.e.~has a section.
A {\bf (proper) solution} to this embedding problem is a finite Galois field extension $F$ of $E$ that contains $F_1$,
together with an isomorphism $\gamma:\Gal(F/E)\rightarrow G$ such that
$\alpha\circ\gamma$ coincides with 
the restriction map from $\Gal(F/E)$ to $\Gal(F_1/E)$.
It was conjectured that over a Hilbertian field $K$ every finite split embedding problem has a solution.

D\`ebes and Deschamps combined these two conjectures into one. 
A {\bf regular solution} to an embedding problem $G\rightarrow\Gal(L/K)$ over $K$
is a solution $\gamma:\Gal(F/K(t))\rightarrow G$ to the induced embedding problem 
$G\rightarrow\Gal(L/K)\cong\Gal(L(t)/K(t))$
over $K(t)$ such that $F/L$ is regular.

\begin{Conjecture}[D\`ebes-Deschamps]\label{DD}
Every finite split embedding problem over any field $K$ has a regular solution.
\end{Conjecture}

This conjecture implies both the regular inverse Galois problem for $K$ (this is the special case $L=K$),
and if $K$ is Hilbertian, then it also implies that every finite split embedding problem over $K$ has a solution.

It is this conjecture that was proven for all ample fields by Pop,
and therefore the reason for the initial interest in ample fields.
Surprisingly, since the publication of Pop's paper in 1996, there was not a single really new case
for which Conjecture \ref{DD} was proven.
In fact, all fields for which this conjecture could be proven, see e.g.~\cite{HarbaterStevenson,Paran},
later turned out to be ample.

In addition, recently doubts came up whether the regular inverse Galois problem 
holds true at all for every field. For example,
in \cite{DebesGhazi} D\`ebes and Ghazi suggest a possible way to disprove the regular inverse Galois problem over $\Q$
(and therefore also to disprove Conjecture \ref{DD}).
Together with the astonishing fact that so far the conjecture could not be proven for any non-ample field,
this led to speculations that 
Conjecture \ref{DD} might actually hold {\em only} for ample fields.
In order to resolve this dispute, we would like to pose the following final question:

\begin{ourque}\label{que:DD}
Does there exist a non-ample field $K$ for which the conjecture of 
D\`ebes and Deschamps holds?
\end{ourque}

A work in this direction is \cite{KoenigsmannIGP},
where a field is constructed for which the regular inverse Galois problem has a positive solution but which does not contain an ample field.
A positive answer to Question~\ref{que:DD} would of course bring up the task to describe 
the class of fields for which Conjecture \ref{DD} holds.

\section*{Acknowledgements}

\noindent
The authors warmly thank 
Elad Paran for helpful comments on the manuscript,
and Jochen Koenigsmann for the permission to publish
his proof connecting Podewski's conjecture to ample fields.

The authors would also like to thank Pierre D\'ebes, Moshe Jarden, Jochen Koenigsmann, Elad Paran, Sebastian Petersen, and Florian Pop
for sharing their thoughts and ideas on ample fields, some of which have found their way into this survey.

The first author is an Alexander von Humboldt postdoc fellow.


\end{document}